\input amstex
\documentstyle{amsppt}
\magnification=\magstep1
\hsize=5.2in
\vsize=6.8in
\centerline {\bf ON THE NOTION OF RELATIVE PROPERTY (T)}
\centerline {\bf FOR INCLUSIONS OF VON NEUMANN ALGEBRAS}
\vskip .1in
\centerline {\rm by}
\vskip .1in
\centerline {\rm JESSE PETERSON and SORIN
POPA\footnote"*"{Supported in part by NSF-Grant 0100883}}

\address Math Dept
UCLA, Los Angeles, CA 90095-155505\endaddress
\email jpete\@math.ucla.edu, \, popa\@math.ucla.edu\endemail

\topmatter

\abstract We prove that the notion of {\it rigidity} (or {\it
relative property} (T)) for inclusions of finite von Neumann
algebras defined in [Po1] is equivalent to a weaker property, in
which no ``continuity constants'' are required. The proof is by
contradiction and uses infinite products of completely positive
maps, regarded as correspondences.
\endabstract

\endtopmatter

\document

The notion of {\it relative property} (T)
(or {\it rigidity}) for inclusions
of finite von Neumann algebras with countable decomposable
center was introduced in ([P1])
by requiring that one of the following conditions (shown equivalent
in [P1]) holds true:

\vskip .05in \noindent {\bf (0.1)}. There exists a normal faithful
tracial state $\tau$ on $N$ such that: $\forall \varepsilon > 0$,
$\exists F'=F'(\varepsilon)\subset N$ finite and $\delta'=
\delta'(\varepsilon) > 0$ such that if $\Cal H$ is a Hilbert
$N$-bimodule with a vector $\xi \in \Cal H$ satisfying the
conditions $\|\langle \cdot\xi, \xi\rangle - \tau\| \leq \delta',
\|\langle \xi \cdot, \xi\rangle - \tau\| \leq \delta'$ and $\|y
\xi - \xi y\| \leq \delta', \forall y\in F'$, then $\exists
\xi_0\in \Cal H$ such that $\|\xi_0 - \xi\| \leq \varepsilon$ and
$b \xi_0 = \xi_0 b, \forall b\in B$. \vskip .05in \noindent {\bf
(0.2)}. There exists a normal faithful tracial state $\tau$ on $N$
such that: $\forall \varepsilon > 0$, $\exists F =
F(\varepsilon)\subset  N$ finite and $\delta =
\delta(\varepsilon)> 0$ such that if $\phi: N \rightarrow N$ is a
normal, completely positive (abreviated c.p. in the sequel) map
with $\tau\circ \phi \leq \tau, \phi(1)\leq 1$ and $\|\phi(x) - x
\|_2 \leq \delta, \forall x\in F$, then $\|\phi(b)-b\|_2 \leq
\varepsilon, \forall b\in B$, $\|b\| \leq 1$. \vskip .05in
\noindent {\bf (0.3)}. Condition $(0.1)$ above is satisfied for
any normal faithful tracial state $\tau$ on $N$. \vskip .05in
\noindent {\bf (0.4)}. Condition $(0.2)$ above is satisfied for
any normal faithful tracial state $\tau$ on $N$. \vskip 0.05in

This definition is the operator algebra analogue of the
Kazhdan-Margulis relative property (T) for inclusions of groups $H
\subset G$ ([M]). It is formulated in the same spirit Connes and
Jones defined the property (T) for single von Neumann algebras in
([CJ]), starting from Kazhdan's property (T) for groups, by using
Hilbert bimodules/c.p. maps, i.e., Connes' {\it correspondences}
([C2]). Thus, while in the case $H = G$ the relative property (T)
of $H \subset G$ amounts to the property (T) of $G$, in the case
$B=N$ and $N$ is a factor the relative property (T) of $B \subset
N$ in the sense of ([P1]) is equivalent to the property (T) of $N$
in the sense of ([CJ]).

But there are in fact two possible ways to define the relative
property (T) for inclusions of groups $H \subset G$: one requiring
that all representations of $G$ that have an almost $G$-invariant
vector must contain a $H$-invariant vector, and another one with
``continuity constants'' of rigidity, requiring in addition that
the vector fixed by $H$ be close to the almost $G$-invariant
vector. The first definition is the original one, formulated in
([M]). The second definition is obtained by adapting to the case
of inclusions of groups a characterization of Kazhdan's property
(T) ``with continuity constants'' which
is implicit in ([DeKi], page 8), ([AW], Lemma 2)
and appears explicitely in ([CJ], Proposition 1) or
([dHV], Proposition 1.16),
a characterization that can be formulated both in terms of unitary
representations and positive definite functions.

Conditions $(0.1)-(0.4)$ are all analogue to this second
definition from group theory. The reason for opting for a
definition ``with continuity constants'' for von Neumann algebras
in ([P1]) is precisely its suitability to a formulation in terms
of completely positive maps (the operator algebra substitute for
positive definite functions), as well as its good behavior to
tensor products and induction/reduction by projections.

For normal subgroups $H \subset G$ the two definitions of relative
property (T) (with and without continuity constants) are easily
seen to be equivalent: the same proof as in the single group case
in (e.g., as in 1.16 of [dHV]) works. But for arbitrary $H \subset
G$, this equivalence is non-trivial and was proved only recently
by Jolissaint ([Jo]). While for applications it is important to
have both definitions available as equivalent conditions, note
that all known examples of subgroups $H \subset G$ with the
relative property (T) are in fact normal (more precisely, $H$
already has the relative property (T) in its normalizer in $G$).

Such equivalence is even more difficult to establish in the
context of von Neumann algebras, where already the single von
Neumann algebra case requires a delicate argument (cf. [CJ]). Yet
it is desirable to have both types of characterizations. Thus,
although for most applications in ([P1]) the characterization
$(0.1)-(0.4)$ is sufficient, a weaker version ``without continuity
constants'' is needed for proving that rigidity is well behaved to
inductive limits (4.5 in [P1]). For this reason, one introduces in
(4.2.2 of [P1]) the notion of $\varepsilon_0$-{\it rigidity},
which requires that \vskip .05in \noindent {\bf (0.5)}. $\exists
F_0 \subset N$ finite and $\delta_0 > 0$ such that if $\phi$ is a
completely positive map on $N$ with $\phi(1) \leq 1$, $\tau \circ
\phi \leq \tau$ and $\|\phi(x)-x\|_2 \leq \delta_0, \forall x \in
F_0$ then $\|\phi(b)-b\|_2 \leq \varepsilon_0, \forall b \in B,
\|b\|\leq 1$, \vskip .1in

\noindent
and one proves in (4.3 of [P1]) that if $N$ is a factor
and $B$ is regular in $N$, i.e., $\Cal N_N(B)\overset \text{\rm
def} \to = \{u \in \Cal U(N) \mid uBu^*=B\}$ generates $N$, then
$B \subset N$ is rigid if and only if it is $1/3${\it-rigid}.

But while enough for proving (4.5 in [P1]),
$\varepsilon_0$-rigidity
is not the exact analogue of the original definition
of relative property (T) for groups ``without
continuity constants'', as considered in ([M]).

The purpose of this paper is to provide such an analogue. Thus, we
prove a characterization of the rigidity for inclusions of finite
von Neumann algebras $B \subset N$ which no longer requires the
invariant vector $\xi_0$ in $(0.1)$ to be close to the almost
invariant vector $\xi$, but merely to be ``almost tracial'' from
left and right (see conditions $(1.2)$, $(1.2')$ in Theorem 1).
This almost traciality condition, which is irrelevant in the group
case, is unavoidable in the framework of von Neumann algebras,
because of the requirement that the rigidity of an inclusion $B
\subset N$ be preserved under reduction by projections in $B$ and
$B'\cap N$ (cf. 4.7 in [P1]).

The almost traciality condition is in fact redundant if we assume
$N$ is a factor, $B'\cap N \subset B$ and $\Cal N_N(B)'\cap N=\Bbb
C$ (see Corollary 2).  This assumption is the same as the one
needed for proving the equivalence between
$\varepsilon_0$-rigidity (condition $(0.5)$) and  rigidity
(conditions $(0.1)-(0.4)$) in (4.3 of [P1]). Consequently, we
obtain a new proof of that result, which is more conceptual and
avoids the use of ultrapower algebras.

\proclaim{Theorem 1} Let $N$ be a finite von Neumann algebra with
countable decomposable center and $B\subset N$ a von Neumann
subalgebra. The following conditions are equivalent: \vskip .05in
\noindent {\bf (1.1)}. The inclusion $B \subset N$ is rigid in the
sense of Definition $4.2$ in ${\text{\rm [Po1]}}$, i.e., it
satisfies the equivalent conditions $(0.1)-(0.4)$. \vskip .05in
\noindent {\bf (1.2)}. There exists a normal faithful tracial
state $\tau$ on $N$ with the property: $\exists F_0 \subset N$
finite and $\delta_0 > 0$ such that $\forall \varepsilon > 0$,
$\exists \delta(\varepsilon)>0$ so that if $\Cal H$ is a Hilbert
$N$-bimodule with a vector $\xi \in \Cal H$ satisfying $\|y \xi -
\xi y\| \leq \delta_0, \forall y\in F_0$, $\|\langle \cdot\xi, \xi
\rangle - \tau\| < \delta, \|\langle \xi \cdot, \xi\rangle -
\tau\| < \delta$, then $\exists \xi_0 \in \Cal H$ such that $b
\xi_0 = \xi_0 b, \forall b\in B$, $\|\langle \cdot\xi_0,
\xi_0\rangle - \tau\| < \varepsilon, \|\langle \xi_0 \cdot,
\xi_0\rangle - \tau\| < \varepsilon$. \vskip .05in \noindent {\bf
(1.2')}. There exists a normal faithful tracial state $\tau$ on
$N$ with the property: $\forall \varepsilon_0 > 0$, $\exists
F_0=F_0(\varepsilon_0) \subset N$ finite and
$\delta_0=\delta_0(\varepsilon_0) > 0$ such that if $\Cal H$ is a
Hilbert $N$-bimodule with a vector $\xi \in \Cal H$ satisfying
$\|y \xi - \xi y\| \leq \delta_0, \forall y\in F_0$, $\|\langle
\cdot\xi, \xi \rangle - \tau\| < \delta_0, \|\langle \xi \cdot,
\xi\rangle - \tau\| < \delta_0$, then $\exists \xi_0 \in \Cal H$
such that $b \xi_0 = \xi_0 b, \forall b\in B$, $\|\langle
\cdot\xi_0, \xi_0\rangle - \tau\| < \varepsilon_0, \|\langle \xi_0
\cdot, \xi_0\rangle - \tau\| < \varepsilon_0$. \vskip .05in
\noindent {\bf (1.3)}. Condition $(1.2)$ above is satisfied for
any normal faithful tracial state $\tau$ on $N$. \vskip .05in
\noindent {\bf (1.3')}. Condition $(1.2')$ above is satisfied for
any normal faithful tracial state $\tau$ on $N$.
\endproclaim

\vskip .1in \proclaim{Corollary 2} Let $N$ be a finite factor and
$B \subset N$ a von Neumann subalgebra such that $B'\cap N \subset
B$, $\Cal N_N(B)'\cap N = \Bbb C$. The following conditions are
equivalent: \vskip .05in \noindent {\bf (2.1)}. The inclusion $B
\subset N$ is rigid. \vskip .05in \noindent {\bf (2.2)}. The
inclusion $B \subset N$ is $\varepsilon_0$-rigid for some $0<
\varepsilon_0 <1$ (i.e., it satisfies condition $(0.5)$ for that
$\varepsilon_0$). \vskip .05in \noindent {\bf (2.3)}. $\exists F_0
\subset N$ finite and $\delta_0 > 0$ such that if $\Cal H$ is a
Hilbert $N$-bimodule with a unit vector $\xi \in \Cal H$
satisfying $\|y \xi - \xi y\| \leq \delta_0, \forall y\in F_0$,
$\|\langle \cdot \xi, \xi \rangle - \tau \| \leq \delta_0$,
$\|\langle \xi, \xi \cdot \rangle - \tau \| \leq \delta_0$, then
$\exists \xi_0 \in \Cal H$, $\xi_0 \neq 0$, such that $b \xi_0 =
\xi_0 b, \forall b\in B$.
\endproclaim

The proof of Theorem 1 proceeds by contradiction, assuming
$(1.2')$ is satisfied while $(1.1)$ is not. By using the point of
view of correspondences (i.e., going back and forth from c.p. maps
to Hilbet bimodules [C2]) and technical background from ([P1]),
this allows us to construct a sequence of completely positive,
subunital, subtracial maps $\phi_n$ on $N$ that get closer and
closer to $id_N$ in point $\|\cdot \|_2$-topology, in a way that
makes the infinite product (composition) $\phi = \Pi_n \phi_n$ be
well defined, close to $id_N$ on prescribed finite subsets of $N$
and with a ``controlled divergence'' from $id_B$, when restricted
to $B$. Moreover, we show that the c.p. maps $\phi_n$ can be taken
so that the operators they induce on $L^2(N, \tau)$ are positive.
If $(\Cal H_\phi, \xi_\phi)$ is the pointed Hilbert $N$-bimodule
associated with $\phi$ as in (1.1 of [P1]), then $\xi_\phi$ almost
commutes with $N$. Thus, by $(1.2')$, $\Cal H_\phi$ contains a
non-zero $B$-central vector $\xi_0$. Approximating $\xi_0$ by
$\Sigma_j x_j \xi_\phi y_j \in {\text{\rm sp}} N\xi_\phi N$, and
using the infinite product form of $\phi$, as well as the
positivity (as operators) of the $\phi_n$'s, leads to a
contradiction.

We mention that Theorem 1 doesn't directly entail the analogous
result for groups in ([Jo]). However, when translating its proof
to the case of inclusions of groups, we obtain a more direct and
shorter proof of ([Jo]), which we present in the Appendix to this
paper. Nevertheless, part of the proof of Theorem 1 was inspired
by ([Jo]). Thus, our idea of using infinite products of c.p. maps
was triggered by an effort to bypass negative definite functions
(for which no satisfactory operator algebra analogue exists) and
their infinite sums, used in ([Jo]). For more on infinite products
of correspondences, both as c.p. maps and as Hilbert bimodules,
see ([Pe]).

\vskip .05in
\noindent
{\it Proof of Theorem 1}.
We clearly have $(1.3') \implies (1.3) \implies(1.2)$
and $(1.3') \implies (1.2') \implies (1.2)$. Also,
condition $(0.3)$
implies $(1.3')$ above, showing that $(1.1)  \implies (1.3')$.
Thus, in order to finish the proof of the theorem it
is sufficient to prove that if
$(1.2)$ is satisfied for a certain normal
faithful tracial state $\tau$ then
condition $(0.2)$ is satisfied for that $\tau$.

To this end, we first
reduce this implication to the case $N$ is separable.
Thus, we begin by proving that $(1.2)$
implies $B$ is separable
(in the norm $\|\quad \|_2$ given by $\tau$).

Let $F = \cup_n F_0(1/n)$ and denote
by $N_0$ the von Neumann algebra generated by $F$, which thus follows separable.
Denote $\Cal H = L^2(\langle N,  N_0 \rangle, Tr)$ and $\xi = e_{N_0}$. Then
$\Cal H$ is a $\langle N, N_0 \rangle$ Hilbert bimodule, so in particular
it is a $N$ bimodule. Since $[N_0, \xi]=0$ and $\xi$
is tracial, by $(1.2)$ it follows that there exists a unit vector $\xi_0 \in \Cal H$ such that
$[B, \xi_0]=0$ and
$Tr(\cdot \xi_0\xi_0^*)$ close to $\tau_N$. It follows that $B$ commutes
with $\xi_0\xi_0^* \in L^1(\langle N, N_0 \rangle, Tr)$, so it also commutes with
the spectral projections of $\xi_0 \xi_0^*$. If $B$ is non-separable, then there exists
$z_0 \in \Cal P(\Cal Z(B))$ such that $Bz$ is non-separable, $\forall z\in \Cal P(\Cal Z(B))$, $z\leq z_0$.
By the condition $Tr(\cdot \xi_0\xi_0^*)$
close to $\tau_N$, it follows that there
exists a spectral projection $e$ of $\xi_0\xi_0^*$ such that $z_0e\neq 0$.
Thus, $z_0e \in \langle N, N_0 \rangle$ commutes with $B$ and has finite trace $Tr$. By
(2.2 in [P2]), it follows that for some $p\in \Cal P(Bz_0)$
and $\eta \in \Cal H$ with $p\eta \neq 0$
we have $pBp\eta \subset \overline{\eta N_0}$.
But this gives a contradiction, since the closure of $\eta N_0$ in $\Cal H$ is a separable
Hilbert space while $pBp \eta$ is not.

Now note that if $(1.2)$ holds true for
$(B \subset N, \tau)$ then it holds true for
$(B \subset N_1, \tau)$, where $N_1$ is the (separable) von Neumann algebra
generated by $B$ and $F$. Indeed, this is immediate to see by inducing
$N_1$-bimodules
to $N$-bimodules. Thus,
if we assume $(1.2) \implies (0.2)$
for inclusions of separable von Neumann algebras, then it follows that
$B \subset N_1$ satisfies the equivalent conditions $(0.1)-(0.4)$.
But then (4.6 in [P1]) shows that $B \subset N$ satisfies
these conditions as well.

Thus, from now on we may assume $N$ is separable. We need
the following equivalent characterization of $(0.2)$.

\proclaim{Lemma 3}  Condition $(0.2)$ for $(B \subset N, \tau)$
holds true if and only if it holds true for completely positive maps $\phi$ with
the property $T_\phi \geq 0$.
\endproclaim

\noindent
{\it Proof}. We first prove that if $(B\subset N, \tau)$ satifies:
\vskip .1in
\noindent
$(3.1)$. $\forall \varepsilon > 0$, $\exists F'
= F'(\varepsilon)\subset  \Cal U(N)$
finite and $\delta' = \delta'(\varepsilon)> 0$ such that if
$\phi: N \rightarrow N$ is a
normal, completely positive map with $\phi^*=\phi$, $\tau\circ \phi \leq \tau,
\phi(1)\leq 1$ and $\|\phi(x) - x \|_2 \leq \delta', \forall x\in
F'$, then $\|\phi(b)-b\|_2 \leq \varepsilon, \forall b\in B$,
$\|b\| \leq 1$.
\vskip .1in
\noindent
then $(B\subset N, \tau), F''(\varepsilon)=F'(\varepsilon^2/2),$
$\delta''(\varepsilon)=\delta'(\varepsilon^2/2)^2/2$ satisfy the condition:
\vskip .1in
\noindent
$(3.2)$. If
$\phi: N \rightarrow N$ is a
normal, completely positive map with $\tau\circ \phi \leq \tau,
\phi(1)\leq 1$ and $\|\phi(x) - x \|_2 \leq \delta'', \forall x\in
F''$, then $\|\phi(b)-b\|_2 \leq \varepsilon, \forall b\in B$,
$\|b\| \leq 1$.
\vskip .1in
Let $\phi$ be as in the hypothesis of $(3.2)$. Then $\phi'=(\phi+\phi^*)/2$ satisfies
$T_{\phi'}^*=T_{\phi'}$ while we still have $\phi'(1)\leq 1$,
$\tau \circ \phi'\leq \tau$. By (1.1.5.3$^\circ$ in [P1]), it follows that if $x\in F''$ then
$$
\|\phi'(x)-x\|_2 \leq
(\|\phi(x)-x\|_2+\|\phi^*(x)-x\|_2)/2
$$
$$
\leq (\|\phi(x)-x\|_2+(2\|\phi(x)-x\|_2)^{1/2})/2
\leq (\delta'' + (2\delta'')^{1/2})/2 \leq \delta'(\varepsilon^2/2)
$$

By $(3.1)$, the above inequality
implies that $\|\phi'(b)-b\|_2 \leq \varepsilon^2/2,$
$\forall b\in B, \|b\| \leq 1$. But then we also have for $b \in B, \|b\|\leq 1$
the estimates:
$$
\|\phi(b)-b\|_2^2 \leq 2\tau(bb^*) - 2 {\text{\rm Re}}\tau(\phi(b)b^*)
$$
$$
= 2\tau(bb^*) - 2 {\text{\rm Re}}\tau(\phi'(b)b^*)\leq 2\|b-\phi'(b)\|_2 \leq \varepsilon^2,
$$
showing that $\|\phi(b)-b\|_2 \leq \varepsilon,$ $\forall b\in B, \|b\|\leq 1$.

Thus, in order to prove Lemma 3 it is now sufficient to show that if
\vskip .1in
\noindent
$(3.0)$. $\forall \varepsilon > 0$, $\exists F_0
= F_0(\varepsilon)\subset  \Cal U(N)$
finite and $\delta_0 = \delta_0(\varepsilon)> 0$ such that if
$\phi: N \rightarrow N$ is a
normal, completely positive map with $T_\phi \geq 0$, $\tau\circ \phi \leq \tau,
\phi(1)\leq 1$ and $\|\phi(x) - x \|_2 \leq \delta_0, \forall x\in
F_0$, then $\|\phi(b)-b\|_2 \leq \varepsilon, \forall b\in B$,
$\|b\| \leq 1$.
\vskip .1in
\noindent
then $(B\subset N, \tau), F_1(\varepsilon)=F_0(\varepsilon/e^2),$
$\delta_1(\varepsilon)=\delta_0(\varepsilon/e^2)/4$ satisfy the condition:
\vskip .1in
\noindent
$(3.1')$. If
$\phi: N \rightarrow N$ is a
normal, completely positive map with $\tau\circ \phi \leq \tau,
\phi(1)\leq 1$, $\phi^*=\phi$
and $\|\phi(x) - x \|_2 \leq \delta_1, \forall x\in
F_1$, then $\|\phi(b)-b\|_2 \leq \varepsilon, \forall b\in B$,
$\|b\| \leq 1$.
\vskip .05in
To show this, let $\phi:N \rightarrow N$ be a c.p. map as in $(3.1')$. Define $\phi''$ on $N$ by
$\phi''=exp(\phi-id_N)= e^{-1} exp(\phi) = e^{-1} \Sigma_{n=0}^\infty \phi^n/n!$,
where $\phi^n$ denotes the $n$ time composition $\phi \circ ... \circ \phi$. By the definition,
it follows that $\phi''$ is completely positive. Also,
since $T_\phi={T_\phi}^*$ and $T_{\phi''} = e^{-1} exp(T_\phi)$, it follows that
$T_{\phi''}$ is a positive operator on the Hilbert space $L^2(N, \tau)$. Moreover, since
$\tau\circ \phi \leq \tau$, we have
$\tau \circ \phi^n \leq \tau, \forall n,$ and thus
$$
\tau\circ \phi'' = e^{-1} \Sigma_n \tau\circ \phi^n/n! \leq \tau\circ \phi \leq \tau.
$$
Similarly, $\phi(1) \leq 1$ implies
$$
\phi''(1) = e^{-1} \Sigma_n \phi^n(1)/n! \leq \phi(1) \leq 1.
$$

By taking into account that
$$
\|\phi''(x)-x\|_2 = \|\Sigma_{n=0}^\infty (\phi-id_N)^n(x)/n! - x \|_2
=\|\Sigma_{n=1}^\infty(\phi-id_N)^n(x)/n!\|_2
$$
$$
\leq (\Sigma_{n=1}^\infty 2^{n-1}/n!) \|\phi(x)-x\|_2
= ((e^2-1)/2) \|\phi(x)-x\|_2 < 4 \|\phi(x)-x\|_2,
$$
it follows that if $\phi$ satisfies the hypothesis of  $(3.1')$ then
$\phi''$ satisfies the conditions in $(3.0)$ for
$\varepsilon$ of the form $e^{-2}\varepsilon$.
Thus, by $(3.0)$ we have
$\|\phi''(b)-b\|_2 \leq e^{-2}\varepsilon,$ $\forall b\in B, \|b\|\leq 1$.
To obtain from this that $\|\phi(b)-b\|_2$ is uniformly small
for $b$ in the unit ball of $B$,
denote $f(t)=1-e^{-t}, 0\leq t\leq 2$.
Since $f(0)=0, f'(t) = e^{-t} \geq e^{-2}$,  it follows that $f(t) \geq te^{-2}, \forall
0\leq t \leq 2$. Hence, if we let $S=1-T_\phi \in \Cal B(L^2(N, \tau))$
then $0 \leq S \leq 2$, which by functional calculus gives $ f(S) \geq e^{-2}S$.
We thus get the estimates:
$$
e^{-2} \varepsilon \geq \|\phi''(b)-b\|_2 = \|f(S)(\hat{b})\|_2
$$
$$
\geq e^{-2} \|S(\hat{b})\|_2
= e^{-2} \|b - \phi(b)\|_2
$$
for all $b$ in $B$ with $\|b\| \leq 1$.

\hfill $\square$

\vskip .05in
\noindent
{\it Proof of} $(1.2) \implies (0.2)$. We proceed by contradiction, assuming that
$(1.2)$ holds true while $(0.2)$ doesn't. Note first that ``non-$(0.2)$''
and the separability of $N$ implies:

\vskip .05in
\noindent
$\overline{(0.2)}$. There exist $c_0 > 0$, c.p. maps
$\{\phi_n\}_n$ on $N$ and unitary
elements $\{b_n\}_n\subset \Cal U(B)$,
such that $T_{\phi_n}\geq 0,$ $\tau \circ \phi_n \leq \tau$,
$\phi_n(1) \leq 1$, $\|\phi_n(x)-x\|_2 \rightarrow 0$,
$\forall x\in N$, and $\|\phi_n(b_n)- b_n\|_2 \geq c_0, \forall n$.

\vskip .05in

Note also that the inequality $\|\phi_n(b_n)- b_n\|_2 \geq c_0$ implies
$$
c_0^2 \leq \|\phi_n(b_n)- b_n\|_2^2 \leq 2 - 2 \tau(\phi_n(b_n)b_n^*)
$$
\noindent
which in turn gives
$$
\|\phi_n(b_n)\|_2^2 =  \tau(\phi_n^2 (b_n)b_n^*)
\leq \tau(\phi_n (b_n)b_n^*) \leq 1-c_0^2/2.
$$

Thus, if we let $c_1 = (1-c_0^2/2)^{1/2}<1$, then $\overline{(0.2)}$
implies:
\vskip .1in
\noindent
$\overline{(3.2)}$. There exist $c_1 < 1$, completely positive maps
$\{\phi_n\}_n$ on $N$ and
$\{b_n\}_n\in \Cal U(B)$,
such that $T_{\phi_n}$ are positive operators on
$L^2(N, \tau)$,  $\tau \circ \phi_n \leq \tau$,
$\phi_n(1) \leq 1$, $\|\phi_n(x)-x\|_2 \rightarrow 0$,
$\forall x\in N$, and $\|\phi_n(b_n)\|_2 \leq c_1, \forall n$.
\vskip .1in

From this point on, we'll need the following notation: If $Y=\{(y_j, z_j)\}_j \subset N \times N$
is a finite set and $\xi\in \Cal H$ for some
$N$ bimodule $\Cal H$ then we denote
$\Cal T_Y(\xi)=
\Sigma_j y_j \xi z_j$. Also, if $\phi:N \rightarrow N$ is a linear map then
$\Cal T_Y \phi(x) = \Sigma_j z_j^*\phi(y_j^* x y_i)z_i$. Note that
if $\phi$ is the normal c.p. map on $N$
given by a vector  ``bounded from the right'' $\xi$
in a Hilbert $N$-bimodule $\Cal H$, as in (1.1.3 of [P1]), then $\Cal T_Y\phi$
is the normal c.p. map given by the vector $\Cal T_Y(\xi)$.

Let $\varepsilon >0$ be so that $\varepsilon < (1-c_1)/(1+c_1)$,
where $c_1$ is as given by $\overline{(3.2)}$. Let $F_0$,
$\delta_0, \delta=\delta(\varepsilon^2/16)$ be as given by $(1.2)$,
with $1 \in F_0$.
Choose a countable $\Bbb Q+ i \Bbb Q$ $\|\quad \|_2$-dense subalgebra
$N_0 =\{x_n\}_n \subset N$.
Denote by $\{Y_n\}_n$ the collection of
all finite sets of pairs of elements in
$N_0$ (i.e., $\{Y_n\}_n$ is the set of finite subsets of $N_0 \times N_0$).
Let $\{\phi_n\}_n$ be the c.p. maps given by
$\overline{(3.2)}$. We choose recursively an increasing sequence of
finite sets $S_k \subset N$ and an increasing sequence of
integers $n_k$, $k \geq 0$,  such that
if we set $S_0=F_0$, $n_0=0, \phi_0=id_N$, $Y_0=\{(1,1)\}$
then for all $k\geq 1$:
$$
\{b_{n_i} \mid 1\leq i \leq k-1\},
\{x_j \mid 1\leq j\leq k\} \subset  S_k \tag a
$$
$$
\Cal T_{Y_j}(S_{k-1}) \subset S_k, \forall j\leq k \tag b
$$
$$
(\phi_{n_{k-1}} \circ ... \circ \phi_{n_j} \circ \Cal T_{Y_i}
(\phi_{n_{j-1}} \circ... \circ \phi_{n_1}))(S_{k-1}) \subset S_k, \forall
j \leq k-1, i \leq k \tag c
$$
$$
\|\phi_{n_k}(x)-x\|_2 \leq \delta^2_0/2^{k+4},
\forall x \in S_k, \|\phi_{n_k}(1)-1\|_2 \leq \delta^2/2^{k+4} \tag d
$$
$$
\|\phi_{n_k}(b_{n_k})\|_2 \leq c_1, \forall k \geq 1 \tag e
$$

\vskip .05in
\noindent

By $(a)-(d)$, it follows that
$\|\phi_{n_k}\circ ... \circ \phi_{n_1}(x_j)-\phi_{n_{k-1}} \circ ... \circ
\phi_{n_1}(x_j)\|_2 \leq \delta_0^2  2^{-k-4}$, $\forall j\leq k$. Thus,
$\{\phi_{n_k}\circ ... \circ \phi_{n_1}(x_j)\}_k$ is Cauchy for each $j$. Since
$\tau\circ \phi_{n_k} \circ ... \circ \phi_{n_0} \leq \tau$ and
$\phi_{n_k} \circ ... \circ \phi_{n_0} (1) \leq 1$, by the density of $\{x_j\}_j$ in $N$
it follows that in fact
$\{\phi_{n_k}\circ ... \circ \phi_{n_1}(x)\}_k$ is Cauchy $\forall x\in N.$
Thus, we can define $\phi(x)=\underset k \rightarrow \infty \to \lim \phi_{n_k}
\circ ... \circ \phi_{n_0}(x)$, $\forall x\in N$, which follows c.p. on $N$ with
$\tau \circ \phi \leq \tau, \phi(1) \leq 1$. Moreover, we have
$$
\|\phi(x) - x\|_2 \leq \Sigma_{k=0}^\infty \|\phi_{n_k} \circ ... \circ \phi_{n_0} (x) -
\phi_{n_{k-1}} \circ ... \circ \phi_{n_0}(x)\|_2
$$
$$
\leq \delta^2_0 \Sigma_{k=1}^\infty 2^{-k-4}  = \delta^2_0/16,
$$
\noindent
for all $x \in F_0$. Similarly, $\|\phi(1)-1\|_2 \leq \delta^2/16$.
Thus, if we denote $(\Cal H_\phi, \xi_\phi)$ the pointed Hilbert bimodule
associated with $\tau(\phi(1))^{-1} \phi$ then by (1.1.2.4$^\circ$
in [P1]) we have
$$
\| x\xi_\phi
-\xi_\phi x\|^2 \leq 2\|\tau(\phi(1))^{-1}\phi(x) - x \|_2^2 +
2 \|\tau(\phi(1))^{-1}\phi(1)\|_2 |\tau(\phi(1))^{-1}\phi(x) - x \|_2
$$
$$
\leq 2 (1-\delta^2/16)^{-2} ((\|\phi(x)-x\|_2 + \|\phi(1)-1\|_2)^2
+ (\|\phi(x)-x\|_2 + \|\phi(1)-1\|_2))
$$
$$
\leq 2 (1-\delta^2/16)^{-2} 2(\delta_0^2/16 + \delta^2/16)
\leq \delta_0^2,
$$
for all $x\in
F_0$, whenever $\delta \leq \delta_0 < 1$.
Also, since $\tau\circ \phi \leq \tau$, by (1.1.5.3$^\circ$ in [P1]) we get
$$
\|\langle \cdot \xi_\phi, \xi_\phi \rangle - \tau\|
= \|\phi^*(1)-\tau(\phi^*(1))\|_1 \tau(\phi(1))^{-1}
$$
$$
2 (1-\delta^2/16)^{-1} \|\phi^*(1) - 1\|_2
\leq 2 (1-\delta^2/16)^{-1} (2 \|\phi(1)-1\|_2)^{1/2}
$$
$$
< 3 (1-\delta^2/16)^{-1} \delta/4 < \delta.
$$

Similarly, we have
$$
\|\langle \xi_\phi \cdot, \xi_\phi \rangle - \tau\| =
\tau(\phi(1))^{-1} \|\phi(1)- \tau(\phi(1))\|_1
$$
$$
\leq 2 (1-\delta^2/16)^{-1} \|\phi(1)-1\|_2 < \delta.
$$
Thus, $(1.2)$ applies to get a unit vector $\xi_0 \in \Cal H_{\phi}$
such that $b\xi_0 = \xi_0 b, \forall b\in B$ and $\|\langle \cdot
\xi_0, \xi_0 \rangle - \tau\| < \varepsilon^2/16,$
$\|\langle \xi_0 \cdot , \xi_0 \rangle - \tau\| < \varepsilon^2/16$.

We'll now use the density of the set $\{x_n\}_n$ in $N$ and
Kaplansky's theorem to show that there exists a vector
of the form $\Cal T_{Y_n} \xi_\phi$ which satisfies the same properties
as $\xi_0$:

\proclaim{Lemma 4} There exists $Y\subset N_0 \times N_0$
such that $\Cal T_Y \phi$ and
$\xi = \Cal T_{Y}(\xi_\phi)$ satisfy $\|b\xi-\xi b\|\leq
\varepsilon,$ $\|\Cal T_Y \phi(b)-b\|_2 < \varepsilon$,
$\forall b\in \Cal U(B)$,
$\|\langle \cdot
\xi, \xi \rangle - \tau\| < \varepsilon,$
$\|\langle \xi \cdot , \xi \rangle - \tau\| < \varepsilon$.
\endproclaim
\vskip .05in
\noindent
{\it Proof}. Since $\{\Cal T_{Y_n}(\xi_\phi)\}_n =
{\text{\rm sp}} N_0 \xi_\phi N_0$ is dense
in sp$N \xi_\phi N$ which  in turn is dense in $\Cal H_\phi$,
it follows that there exists $n_0$ such that
if we denote $Y'=Y_{n_0}$, then
$\xi = \Cal T_{Y'}(\xi_\phi)$ is close enough
to $\xi_0$ to ensure that $\|b\xi-\xi b\|\leq
\varepsilon^2/16, \forall b\in \Cal U(B)$, while we still have
$\|\langle \cdot
\xi, \xi \rangle - \tau\| < \varepsilon^2/16,$
$\|\langle \xi \cdot , \xi \rangle - \tau\| < \varepsilon^2/16$.

Let $a_{0} = \Cal T_{Y} \phi(1)$,
$d_{0} = (\Cal T_{Y} \phi)^*(1)$. Then the last two
inequalities become $\|a_{0} -1 \|_1 < \varepsilon^2/16$,
$\|d_{0} -1 \|_1 < \varepsilon^2/16$.
Since $\xi$ implements the c.p. map
$\Cal T_{Y'}(\phi)$, by (Lemma 1.1.3 in [P1]) it follows that if we let
$a = a_{0}\vee 1$, $d = d_{0}\vee 1$ and $\xi'= a^{-1/2} \xi d^{-1/2}$
then $\|\xi-\xi'\|^2 \leq 8 \varepsilon^2/16=\varepsilon^2/2$.

By Kaplansky's density theorem
there exist $a'_n, d'_n $ in the unit ball of
$N_0$ such that
$\underset n \rightarrow \infty \to \lim \|a'_n - a^{-1/2}\|_2 = 0$,
$\underset n \rightarrow \infty \to \lim \|d'_n - d^{-1/2}\|_2 = 0$.

Denote $Y'_n = \{(a'_n x, y d'_n) \mid (x,y) \in Y'\}$,
$\xi'_n = \Cal T_{Y'_n} (\xi_\phi) = a'_n \xi d'_n$,
$\phi' = \phi_{\xi'}$,
$\phi'_n = \Cal T_{Y'_n} \phi$ and note that
$Y'_n \in N_0\times N_0$, $\forall n$. It follows that
$\underset n \rightarrow \infty \to \lim \|\xi'-\xi'_n\|=0$.
Also, it is easy to see that $\|\phi_n'(1)\| \leq \|a_0\|$.
This implies
$\|\phi_n'(1) - \phi'(1)\|_2 \rightarrow 0$,
so in particular $\|\phi_n'(1)\|_2
\rightarrow  \|\phi'(1)\|_2$.  But  by (1.1.3 in [P1])
we have
$$
\|\phi'_n(b)-b\|^2_2 \leq \|[b, \xi'_n]\|^2 + (\|\phi'_n(1)\|_2^2 -1)
$$
and since
$$
\|[b, \xi'_n]\| \leq \|[b, \xi]\|+\|\xi-\xi'\|+\|\xi' - \xi'_n\|
$$
$$
\leq \varepsilon^2/16 + 2^{-1/2} \varepsilon + \|\xi' - \xi'_n\|,
$$
for large enough $n$ and $\varepsilon < 1$ (to insure that
$\varepsilon^2/16 + 2^{-1/2} \varepsilon < \varepsilon$)
we obtain the estimate
$\|\phi'_n(b)-b\|_2 <  \varepsilon, \forall b\in \Cal U(B)$.

Thus, if we put $Y=Y'_n$ then all the requirements are satisfied.

\hfill $\square$

For each $1\leq j \leq k$, denote
$\phi_j^k = \phi_{n_k} \circ ... \circ \phi_{n_j}$,
$\phi_j^\infty = \underset k \rightarrow \infty \to \lim \phi_j^k$.
Since $\phi_{n_k} \rightarrow id_N$, by (Corollary 1.1.2 in [P1])
it follows that $\|\phi_{n_k}(y^* \cdot z) - y^* \phi_{n_k}( \cdot) z \|
\rightarrow 0$, $\forall y,z\in N$, and thus
$$
\underset k \rightarrow \infty \to \lim \| \Cal T_Y \phi
-  \phi_{k+1}^\infty (\Cal T_Y \phi_1^k) \| =0 \tag f
$$

Since by $(a), (b), (c)$
we have $\phi_{m+1}^j (\phi_{k+1}^m(\Cal T_Y \phi_1^k(b_{n_m}))) \in S_j$,
$\forall j > m$, we also get
$$
\underset m \rightarrow \infty \to \lim
\|\phi_{k+1}^\infty (\Cal T_Y \phi_1^k)(b_{n_m})-
\phi_{k+1}^m (\Cal T_Y \phi_1^k)(b_{n_m})\|_2=0 \tag g
$$

Altogether, from $(f)$ and $(g)$ we obtain:
$$
\|b_{n_m} - \phi_{k+1}^m (\Cal T_Y \phi_1^k)(b_{n_m})\|_2
$$
$$
\leq \|b_{n_m} - \phi_{k+1}^\infty (\Cal T_Y \phi_1^k)(b_{n_m})\|_2
+ \|\phi_{k+1}^\infty (\Cal T_Y \phi_1^k)(b_{n_m})-
\phi_{k+1}^m (\Cal T_Y \phi_1^k)(b_{n_m})\|_2
$$
$$
\leq \|b_{n_m} -\Cal T_Y\phi(b_{n_m})\|_2 + \|\Cal T_Y \phi
-  \phi_{k+1}^\infty (\Cal T_Y \phi_1^k) \| +
\varepsilon (m)
$$
$$
\leq \|b_{n_m} -\Cal T_Y\phi(b_{n_m})\|_2 + \varepsilon' (k) + \varepsilon (m)
$$
$$
< \varepsilon + \varepsilon' (k) + \varepsilon (m)
$$
for all $m > k$,
with $\varepsilon(m), \varepsilon'(k)$ satisfying
$\underset m \rightarrow \infty \to \lim \varepsilon(m)= 0$,
$\underset k \rightarrow \infty \to \lim \varepsilon'(k)=0$.
Thus, there exists $k_0$ such that for all $m > k \geq k_0$ we have
$$
\|b_{n_m} - \phi_{k+1}^m (\Cal T_Y \phi_1^k)(b_{n_m})\|_2
\leq \varepsilon \tag h
$$

Moreover, since by Lemma 4 we have
$\|\Cal T_Y \phi(1)\|_2 < 1+ \varepsilon$,
by $(f)$ it follows that we can choose $k_0$ so that
for all $k \geq k_0$ we also have
$$
\|\phi_{k+1}^\infty (\Cal T_Y \phi_1^k)(1)\|_2
< 1 + \varepsilon  \tag i
$$

Fix some $k \geq k_0$. By $(g)$ and $(i)$
it follows that there exists $m > k$
such that
$$
\|\phi_{k+1}^{m-1} (\Cal T_Y \phi_1^k)(1)\|_2 < 1+\varepsilon \tag j
$$

On the other hand, since $b_{n_m}$ are unitary elements, by (part
$1^\circ$ in Lemma 1.1.2 of [P1]) we have
$$
\|\phi_{k+1}^{m-1} (\Cal T_Y \phi_1^k)(b_{n_m})\|_2 \leq
\|\phi_{k+1}^{m-1} (\Cal T_Y \phi_1^k)(1)\|_2 \tag k
$$

Combining $(j)$ and $(k)$ we get
$$
\|\phi_{k+1}^{m-1} (\Cal T_Y \phi_1^k)(b_{n_m})\|_2  \leq 1+ \varepsilon \tag l
$$

For simplicity, denote by $T$ the operator implemented by
$\phi_{k+1}^{m-1} (\Cal T_Y \phi_1^k)$ on $L^2(N, \tau)$,
by $S$ the operator implemented by $\phi_{n_m}=\phi_{m}^m$ on
$L^2(N, \tau)$ and by $\eta$ the vector $\hat{b}_{n_m}$. Recall that
$S\geq 0$ (so in particular $S=S^*$). By $(e)$ we have $\|S(\eta)\|_2
\leq c_1$ and by $(l)$ we have $\|T(\eta)\|_2 \leq 1+ \varepsilon$.
Also, by $(h)$ we have $\|\eta-ST(\eta)\|_2 \leq \varepsilon$.
By applying
twice the Cauchy-Schwartz inequality, it follows that
$$
\varepsilon \geq | \langle \eta - ST(\eta), \eta \rangle |
$$
$$
= | 1 - \langle T(\eta), S (\eta) \rangle |
\geq 1 - \|T(\eta)\|_2 \|S (\eta)\|_2
\geq 1 - c_1 (1+\varepsilon)
$$

But this implies $\varepsilon \geq (1-c_1)/(1+c_1)$, in
contradiction with our  initial choice of $\varepsilon$.

\hfill $\square$

\vskip .1in
\noindent
{\it Proof of Corollary 2}. We clearly have $(2.1) \implies (2.2)$.
To prove $(2.2)
\implies (2.3)$, let $(\Cal H, \xi)$ be a pointed Hilbert $N$-bimodule
satisfying $(2.3)$ for some $F_0 \subset N$ finite $\delta_0>0$. Denote
$a_0, b_0\in L^1(N, \tau)$ the Radon-Nykodim derivatives
of $\langle \cdot \xi, \xi \rangle$ and $\langle \xi \cdot, \xi \rangle$
with respect to $\tau$. Let $\xi' = (a_0 \vee 1){-1/2} \xi
(b_0 \vee 1)^{-1/2}$
and note that $\xi'$ implements subtracial functionals on $N$,
both left and right. By
(1.1.3.1$^\circ$ in [P1]) we have $\|\xi - \xi'\|^2 \leq
8 \delta_0$. Also, since
$\|(a_0 \vee 1)^{-1}\xi\|\leq 1$, $\xi (b_0 \vee 1)^{-1}\|\leq 1$,
by applying the Cauchy-Schwartz inequality we get
$$
\|\xi'\|^2=\langle
(a_0 \vee 1)^{-1}\xi, \xi (b_0 \vee 1)^{-1} \rangle
$$
$$
\geq \langle \xi, \xi  \rangle
- \|(a_0 \vee 1)^{-1}\xi - \xi\| - \|\xi - \xi (b_0 \vee 1)^{-1}\|
$$
$$
\geq 1 - (\tau((1-(a_0 \vee 1)^{-1})^2) + \|a_0 -1 \|_1)^{1/2}
-  (\tau((1-(b_0 \vee 1)^{-1})^2) + \|b_0 -1 \|_1)^{1/2}
$$
$$
\geq 1 - (2 \|a_0 -1 \|_1)^{1/2} + (2 \|b_0 -1 \|_1)^{1/2}
\geq 1- 3 \delta_0^{1/2}.
$$

Denote $\phi'=\phi_{(\Cal H,  \xi'')}$
the c.p. map associated with $\xi''=\|\xi'\|^{-1} \xi'$
as in (1.1.3 in [P1]). By
(1.1.3.1$^\circ$ in [P1]) we have $\|\xi - \xi''\|^2 \leq
8 (1-3\delta_0^{1/2})^{-1} \delta_0$. Also, by
(1.1.3.2$^\circ$ in [P1]), if we take $\delta_0 \leq 1/16$
then we have:
$$
\|\phi'(u) - u\|_2
\leq \|[u, \xi'']\|^2 + (\|\phi'(1)\|_2^2-1)
$$
$$
\leq \|[u, \xi']\|^2  + (1-3\delta_0^{1/2})^{-2}-1
$$
$$
\leq (\|[u, \xi]\| + \|\xi-\xi'\|)^2 + 100 \delta_0^{1/2},
$$
for all unitary elements in $N$, and thus for all
$x$ in the unit ball of $N$. Thus, if we take
$\delta_0$ sufficiently small and apply
the inequality to $x \in F_0$,
then $\phi'$ checks condition $(2.2)$. Thus,
$\|\phi'(b)-b\|_2
\leq \varepsilon_0, \forall b\in \Cal U(B)$.
But then we have the estimates
$$
\varepsilon_0 \geq |\langle b-\phi'(b), b \rangle | = | 1 -
\tau (\phi'(b)b^*)|
$$
$$
=|\langle \xi'', (\xi''  -  b\xi'' b^*) \rangle |.
$$

Thus, if we let $\xi_0$ be element of minimal norm
in $\overline{\text{\rm co}}
\{ b \xi'' b^* \mid b \in \Cal U(B) \} \subset \Cal H$
then $\xi_0 \neq 0$ and $[\xi_0, B]=0$. This ends the proof
of $(2.2) \implies (2.3)$.

By
Theorem 1, in order to prove $(2.3) \implies (2.1)$ it
is sufficient to show that
if $\Cal H$ is a Hilbert $N$-bimodule with a
non-zero $B$-central vector $\xi_0 \in \Cal H$, then $\Cal H$ contains a
non-zero $B$-central vector $\xi_1 \in \Cal H$ which is
left and right almost
tracial on $N$. To see this, let $a_0, b_0 \in L^1(N, \tau)_+$
be so that $\tau(xa_0)=\langle x\xi_0, \xi_0 \rangle, \forall x\in N$
and $\tau(b_0x)=\langle \xi_0 x, \xi_0 \rangle, \forall x\in N$.
Since $[B, \xi_0]=0$ and $B'\cap N=\Cal Z(B)$, both $a_0, b_0$
belong to $L^1(\Cal Z(B), \tau)_+$ and have the same support projection.

Thus, given any $\varepsilon > 0$ there exists a non-zero projection
$q \in \Cal Z(B)$ such that $\| a_0q-cq\| < \varepsilon$,
$\| b_0q-cq\| < \varepsilon$, where $c = \langle q
\xi_0, \xi_0 \rangle/\tau(q)$ $=\tau(a_0q)/tau(q)=\tau(b_0q)/\tau(q)$.
Moreover, since $\Cal N(B)$ acts ergodically on $\Cal Z(B)$ (because
$\Cal N(B)'\cap '\cap N =\Bbb C$), it follows that we can take
$q$ to satisfy $\tau(q)=1/n$, for some integer $n\geq 1.$ Furthermore,
there exist partial isometries $v_j \in \Cal G\Cal N(B)$, $1\leq j\leq n$,
such that $v_j^*v_j = q, \forall j$, $\Sigma_j v_jv_j^* = 1$.
But then $\xi_1 = c^{-1/2} \Sigma_j v_j \xi_0 v_j^*$ is easily seen
to satisfy $[\xi_1, B]=0$, $\| \langle \cdot \xi_1, \xi_1 \rangle
-\tau\| \leq \varepsilon$,  $\| \langle \xi_1 \cdot, \xi_1 \rangle
-\tau\| \leq \varepsilon$.

\hfill $\square$

\noindent
{\bf Remarks 5}. $1^\circ$. The equivalence in Corollary 2
can actually be proved without assuming $B'\cap N \subset B$,
but the argument becomes considerably longer.

2$^\circ$. Another proof of Lemma 3 can be obtained by
following the argument on page (40 of [P1]), which shows that
if $\{\phi_n\}_n$ are the c.p. maps given by $\overline{(0.2)}$
then for ``large'' $k$ and ``very large'' $n$ the average map
$\psi=k^{-1} \Sigma_{j=1}^k \phi_n^j$ shrinks the norm $\|\cdot \|_2$
of some elements $b$ in the unit ball of $B$ by a
uniform constant $c_1 < 1$. But then $\phi=\psi^* \circ \psi$
still shrinks $b$ by $c_1$
(thus also some unitaries in $B$, by the
Russo-Dye Theorem) and is positive as an operator on $L^2(N, \tau)$.

3$^\circ$. By (4.7.2$^\circ$ in [P1]), if an inclusion of finite
von Neumann algebras with countable decomposable center $B \subset
N$ satisfies any of the equivalent conditions in Theorem 1 then it
satisfies the following {\it uniform local weak rigidity}
condition (u.l.w.r.): ``For all projections $p\in \Cal P(B) \cup
\Cal P(B'\cap N)$ the inclusion $pBp \subset pNp$ satifies
condition $(2.3)$.'' It would be interesting to prove that,
conversely, if $B \subset N$ satisfies the u.l.w.r. condition then
it is rigid.

\head  Appendix\endhead

We present here a short proof of the result in ([Jo]) showing that
the original definition of the relative property (T) for
inclusions of groups $H \subset G$ in ([M]) is equivalent to a
condition ``with continuity constants'', which can be formulated
both in terms of representations and positive definite functions.
Our proof uses infinite products of positive definite functions,
and is obtained by translating the proof of Theorem 1 to the case
of groups, where many simplifications occur. For simplicity, we
present the proof for discrete groups only, noting that the same
argument works for locally compact groups.

\proclaim{Theorem [Jo]} Let $G$ be a countable discrete group and
$H\subset G$ a subgroup. The following conditions are equivalent:
\vskip .05in \noindent $(A.1)$. $\forall \varepsilon > 0$,
$\exists F=F(\varepsilon) \subset G$ finite and
$\delta=\delta(\varepsilon)>0$ so that if $\varphi:G \rightarrow
\Bbb C$ is a positive definite function and $|\varphi(g) - 1| \leq
\delta, \forall g \in F$ then $|\varphi(h) - 1| \leq \varepsilon,
\forall h \in H$. \vskip .05in \noindent $(A.2)$. $\forall
\varepsilon > 0$, $\exists F'=F'(\varepsilon) \subset G$ finite
and $\delta'=\delta'(\varepsilon)>0$ so that if $\pi:G \rightarrow
\Cal U (\Cal H)$ is a unitary representation of $G$ with a unit
vector $\xi \in \Cal H$ satisfying $\|\pi(g)\xi - \xi\| \leq
\delta', \forall g\in F'$ then $\exists \xi_0 \in \Cal H$ such
that $\pi(h)\xi_0 = \xi_0, \forall h \in H$ and $\|\xi_0 - \xi\|
\leq \varepsilon$. \vskip .05in \noindent $(A.3)$. $\exists F_0
\subset G$ finite and $\delta_0>0$ so that if $\pi:G \rightarrow
\Cal U (\Cal H)$ is a unitary representation of $G$ with a unit
vector $\xi \in \Cal H$ satisfying $\|\pi(g)\xi - \xi\| \leq
\delta_0, \forall g\in F_0$ then $\exists \xi_0 \in \Cal H$,
$\|\xi_0\| = 1$ such that $\pi(h)\xi_0 = \xi_0, \forall h \in H$.
\endproclaim

\vskip .05in
\noindent
{\it Proof.}$(A.1)\implies(A.2)$ is easily seen using the
GNS construction and $(A.2)\implies(A.3)$ is trivial hence
we will proceed with  $(A.3)\implies(A.1)$.

Assume that $(A.3)$ holds but $(A.1)$ does not. Since $(A.1)$ does
not hold we have \vskip .1in \noindent $\overline{(A.1)}$.
$\exists c_0 > 0$, $\{b_n\}_n$ a sequence in $G$, and $\varphi_n:G
\rightarrow \Bbb C$ a sequence of positive definite functions such
that $|\varphi_n(g) - 1| \rightarrow 0, \forall g \in G$, and
$|\varphi(b_n) - 1| \geq c_0, \forall n$. \vskip .1in Moreover, as
$\varphi_n(1) \rightarrow 1$, by substituting
$\varphi_n(1)^{-1}\varphi_n(g)$ for $\varphi_n$ we may assume
$\varphi_n(1) = 1, \forall n$. Also, since $\varphi_n'(g) =
exp(Re(\varphi_n(g)) - 1)$ satisfies $|\varphi_n'(g) - 1|
\rightarrow 0, \forall g \in G$ and $1 - \varphi_n'(b_n) \geq
e^{-2}(1 - Re(\varphi_n(b_n))) \geq e^{-2}/2|\varphi_n(b_n) - 1|^2
\geq e^{-2}/2c_0^2$ (the latter due to the inequality $1 - e^{-t}
\geq e^{-2}t, \forall 0 \leq t \leq 2$), by setting $c_1 =
e^{-2}/2c_0^2$ and substituting $\varphi_n'$ for $\varphi_n$ it
follows that $\overline{(A.1)}$ implies: \vskip .1in \noindent
$\overline{(A.1)'}$. $\exists  c_1 > 0$, $\{b_n\}_n$ a sequence in
$G$, and $\varphi_n:G \rightarrow \Bbb C$ a sequence of positive
definite functions such that $|\varphi_n(g) - 1| \rightarrow 0,
\forall g \in G$, $\varphi_n(1)=1, \varphi_n \geq 0$ and
$|\varphi(b_n) - 1| \geq c_1, \forall n$. \vskip .1in Note that
since $\varphi_n(g) \leq \varphi_n(1) = 1$ we have $\varphi_n(b_n)
\leq 1-c_1, \forall n$. Also if $\{K_n\}_n$ is an increasing
sequence of finite sets in $G$ such that $\cup K_n = G$ then since
$(\varphi_n(b_n))^k \leq (1-c_1)^k, \forall n$ we may construct a
new sequence of positive definite functions by letting
$\varphi_k'=(\varphi_{n_k})^k, \forall k,$ where $n_k >> n_{k-1}$
are chosen increasing rapidly enough to ensure that
$|(\varphi_{n_k}(g))^k - 1| \leq 1/k, \forall g \in K_k$. We then
have that $\varphi_k'$ still satisfy $\overline{(A.1')}$ but also
if $b_k' = b_{n_k}$ then $\varphi_k'(b_k') \rightarrow 0$.

If $Y=\{(\alpha_j, g_j)\}_j
\subset \Bbb C \times G$ is a finite set and $\pi:G \rightarrow
\Cal U (\Cal H)$ is a unitary representation with $\xi\in \Cal H$
then we denote $T_Y(\xi) \overset \text{\rm def} \to =
\Sigma_j \alpha_j \pi(g_j)\xi$. Also, if
$\psi:G \rightarrow \Bbb C$ is a function then $T_Y \varphi(g)
\overset \text{\rm def} \to =
\Sigma_{i,j} \overline{\alpha_j}\alpha_i\varphi(g_j^{-1}gg_i)$.
Note that if $\psi$ is the positive definite function on $G$ given
by a vector $\xi$ in $\Cal H$, then $T_Y\psi$ is the positive
definite function given by the vector $T_Y(\xi)$. Also note that
$T_Y(\psi_1\psi_2)(g) = \Sigma_{i,j}\overline {\alpha_j}\alpha_i
\psi_1(g_j^{-1}gg_i)\psi_2(g_j^{-1}gg_i)$, which is different from
$\psi_1T_Y(\psi_2)$ and $T_Y(\psi_1)T_Y(\psi_2)$.

Let $\varepsilon >0$. Let $F_0$,
$\delta_0$ be as given by $(A.3)$.
Let $\{K_k\}_k$ be an increasing sequence of finite sets in $G$
such that $F_0 \subset K_1$, and $\cup_kK_k = G$.
Let $\{\varphi_k'\}_k$ be the c.p. maps given above. We choose a subsequence
$\{\varphi_{k_j}'\}_j$ of $\{\varphi_k'\}_k$ such that:
$$
|\varphi_{k_j}'(g) - 1| \leq \delta_0^2/2^{j+1},
\forall g \in K_j \tag a
$$
\vskip .05in
\noindent

By $(a)$, it follows that
$|\Pi_{i=1}^j\varphi_{k_i}'(g) - \Pi_{i=1}^{j-1}\varphi_{k_i}'(g)|
\leq \delta_0^2/2^{j+1}, \forall g \in K_j$.  Thus
$\{\Pi_{i=1}^j\varphi_{k_i}'(g)\}_j$ is Cauchy, $\forall g \in G$.
For each $1 \leq n \leq m$, denote
$\varphi_n^m = \Pi_{i=n}^m\varphi_{k_i}'$,
$\varphi_n^\infty = \Pi_{i=n}^\infty\varphi_{k_i}'$,
$\varphi = \Pi_{i=1}^\infty\varphi_{k_i}'$.
Then $\varphi$ is a positive definite function on G such that
$\varphi \geq 0$, $\varphi(1) = 1$, and
$$
|\varphi(g) - 1| \leq \Sigma_{m=1}^\infty |\varphi_1^m - \varphi_1^{m-1}|
\leq \delta_0^2\Sigma_{m=1}^\infty 2^{-m-1} = \delta_0^2/2, \tag b
$$
$\forall g \in F_0$.  Hence if we denote
$(\pi_\varphi:G \rightarrow \Cal U(\Cal H_\varphi), \xi_\varphi)$
the pointed unitary representation associated with $\varphi$, then
$\|\pi_\varphi(g)\xi_\varphi - \xi_\varphi\|^2 = 2 - 2\varphi(g) \leq \delta_0^2$
for all $g \in F_0$.  Thus (A.3) applies to get a unit vector
$\xi_0 \in \Cal H_\varphi$ such that \
$\pi_\varphi(h)\xi_0 = \xi_0, \forall h \in H$.

Since $\xi_\varphi$ is a cyclic vector for $\Cal H_\varphi$ there exists
$Y \subset \Bbb C \times G$ finite such that
$\xi = T_Y(\xi_\varphi)$ satisfies
$\|\xi_0 - \xi\| \leq \varepsilon/3$ and $\|\xi\|=1$.
Then we have
$$
|1 - T_Y(\varphi)(h)| =
|1 - \langle \pi_\varphi(h)\xi,\xi \rangle| \leq
$$
$$
|\langle \pi_\varphi(h)\xi_0,\xi_0 \rangle -
\langle \pi_\varphi(h)\xi_0,\xi \rangle| +
|\langle \pi_\varphi(h)\xi_0,\xi \rangle -
\langle \pi_\varphi(h)\xi,\xi \rangle| \leq \tag c
$$
$$
\| \xi_0 - \xi \| +  \| \xi_0 - \xi \| \leq \varepsilon,
$$
for all $h \in H$.

Using the inequality
$|\varphi'(g_1) - \varphi'(g_2)|^2 \leq 2\varphi'(1)(\varphi'(1) - Re\varphi'(g_1^{-1}g_2))$,
for positive definite functions
$\varphi'$ on $G$ and $g_1,g_2 \in G$, together
with the fact that $\varphi_n^\infty \rightarrow 1$, it follows that
$\|\varphi_n^\infty(g_1 \cdot g_2) - \varphi_n^\infty( \cdot ) \|_\infty
\rightarrow 0, \forall g_1,g_2 \in G$. Thus
$$
\lim_{n \rightarrow \infty} \|T_Y(\varphi) -
\varphi_{n+1}^\infty T_Y(\varphi_1^n)\|_\infty = 0 \tag d
$$
Also since $\varphi_n^\infty(b_{k_m}') \leq \varphi_{k_m}'(b_{k_m}'), \forall m \geq n$
it follows that
$$
\lim_{m \rightarrow \infty} |(\varphi_{n+1}^\infty T_Y(\varphi_1^n))(b_{k_m}')| =
\lim_{m \rightarrow \infty} |\varphi_{n+1}^\infty(b_{k_m}')T_Y(\varphi_1^n)(b_{k_m}')| =  0, \tag e
$$
for all $n \geq 1$.  Combining $(d)$ and $(e)$ we have
$$
\lim_{m \rightarrow \infty} |T_Y(\varphi)(b_{k_m}')| = 0 \tag f
$$
But this contradicts $(c)$ for $\varepsilon < 1$.

\hfill $\square$

\head  References\endhead

\item{[AW]} C. Akeman, M. Walters: {\it Unbounded negative
definite functions} Can. J. Math., {\bf 33} (1981), 862-871.

\item{[C1]} A. Connes: {\it A type II$_1$
factor with countable fundamental group}, J. Operator
Theory {\bf 4} (1980), 151-153.

\item{[C2]} A. Connes: {\it Classification des facteurs},
Proc. Symp. Pure Math.
{\bf 38}
(Amer. Math. Soc., 1982), 43-109.

\item{[CJ]} A. Connes, V.F.R. Jones: {\it Property} (T)
{\it for von Neumann algebras}, Bull. London Math. Soc. {\bf 17} (1985),
57-62.

\item{[DeKi]} C. Delaroche, Kirilov: {\it Sur les relations entre
l'espace dual d'un groupe et la structure de ses sous-groupes fermes},
Se. Bourbaki, 20'e ann\'ee, 1967-1968, no. 343, juin 1968.

\item{[dHV]} P. de la Harpe, A. Valette: ``La propri\'et\'e T
de Kazhdan pour les
groupes localement compacts'', Ast\'erisque {\bf 175}, Soc. Math. de France (1989).

\item{[Jo]} P. Jolissaint: {\it On the relative property T},
preprint 2001.

\item{[K]} D. Kazhdan: {\it Connection of the dual space of a group
with the structure of its closed subgroups}, Funct. Anal. and its Appl.
{\bf1} (1967), 63-65.

\item{[M]} G. Margulis: {\it Finitely-additive invariant measures
on Euclidian spaces}, Ergodic. Th. and Dynam. Sys. {\bf 2} (1982),
383-396.

\item{[Pe]} J. Peterson: {\it On infinite products of
correspondences}, Thesis, UCLA, in preparation.

\item{[P1]} S. Popa: {\it On a class of type} II$_1$ {\it factors with Betti numbers invariants},
preprint OA/0209130.

\item{[P2]} S. Popa: {\it Strong rigidity of} II$_1$
{\it factors coming from malleable actions of weakly rigid groups}, I,
math.OA/0305306.

\item{[Po3]} S. Popa: {\it Correspondences}, INCREST preprint 1986,
unpublished.

\enddocument